\documentclass{article}
\usepackage{amsmath,enumitem,cleveref,amsfonts,amssymb}
\usepackage[pdftex]{graphicx}
\usepackage{caption}
\usepackage{subcaption}
\usepackage{float}
\restylefloat{table}
\restylefloat{figure}

\newtheorem{fact}{Fact}[section]
\usepackage{url}
\usepackage{cases}
\floatstyle{plaintop}
\usepackage[section]{placeins}
\usepackage{hhline}
\usepackage{rotating}
\usepackage{longtable}

\usepackage{array}
\newcolumntype{L}[1]{>{\raggedright\let\newline\\\arraybackslash\hspace{0pt}}m{#1}}
\newcolumntype{C}[1]{>{\centering\let\newline\\\arraybackslash\hspace{0pt}}m{#1}}
\newcolumntype{R}[1]{>{\raggedleft\let\newline\\\arraybackslash\hspace{0pt}}m{#1}}

\newcommand{\T}{\mathcal{T}} 
\newcommand{\B}{\mathcal{B}} 
\newcommand{\W}{\mathcal{W}} 
\newcommand{\G}{\mathcal{G}} 

\renewcommand{\d}{\delta} 
\renewcommand{\S}{\mathcal{S}} 

\graphicspath{{Graphics/}}

\newcommand{\norm}[1]{\lVert #1 \rVert}

\newcommand{\rom}[1]{\uppercase\expandafter{\romannumeral #1\relax}} 
\newcommand{\rpm}{\sbox0{$1$}\sbox2{$\scriptstyle\pm$}
  \raise\dimexpr(\ht0-\ht2)/2\relax\box2 } 

\title{Optimizing horizontal alignment of roads in a specified corridor}
\author{Sukanto Mondal, Yves Lucet\footnote{Computer Science, Arts \& Sciences, 3187 University Way, UBC Okanagan, V1V 1V7, yves.lucet@ubc.ca}, Warren Hare}

\date{\today}

\begin {document}
\maketitle

\begin{abstract}
Finding an optimal alignment connecting two end-points in a specified corridor is a complex problem that requires solving three interrelated sub-problems, namely the horizontal alignment, vertical alignment and earthwork optimization problems. In this research, we developed a novel bi-level optimization model combining those three problems. In the outer level of the model, we optimize the horizontal alignment and in the inner level of the model a vertical alignment optimization problem considering earthwork allocation is solved for a fixed horizontal alignment. Derivative-free optimization algorithms are used to solve the outer problem. The result of our model gives an optimal horizontal alignment in the form of a linear-circular curve and an optimal vertical alignment in the form of a quadratic spline. Our model is tested on real-life data. The numerical results show that our approach improves the road alignment designed by civil engineers by 27\% on average, resulting in potentially millions of dollars of savings.
\end{abstract}


\section{Introduction}
Road design optimization is the problem of finding a curve that minimizes the construction cost while satisfying all of the desired design specifications. The problem is usually divided into three interrelated sub-problems \cite{AASHTO-04}: horizontal alignment optimization, vertical alignment optimization, and earthwork optimization. In typical road design, first a horizontal alignment is proposed and then an associated vertical alignment is optimized considering the earthwork allocations and vertical alignment design constraints.

For a fixed horizontal alignment, the optimization of vertical alignment is a well studied problem \cite{Hare-11,CHENG-13,LIMA-13,Burdett-2014,HARE-14b}.  Many of the past approaches have modeled the optimization of vertical alignment by way of mixed integer linear programming.  This creates a complicated, but deterministic optimization problem that is generally solvable using modern MILP solvers (assuming reasonable road lengths and time allowances) \cite{HARE-14b}.  This implies that given a proposed horizontal alignment, it is possible to evaluate the quality of that alignment in terms of the optimal cost vertical alignment.  Which further suggests that it should be possible to have a computer evaluate and seek an optimal horizontal alignment (in terms of the minimal vertical alignment construction cost).  In this paper, we demonstrate the practicality of this idea, and further demonstrate the value of this approach in terms of cost savings for the final road design.

To do this, we formulate the horizontal alignment optimization problem as a bi-level optimization problem. In the inner level, for a fixed horizontal, a vertical alignment optimization problem is solved using the mixed integer linear programming (MILP) from \cite{Hare-2014} (see Appendix \ref{appendix}), which builds on \cite{HARE-11a,Hare-2014,HARE-14b,Hare-11} and provides a global optimum. The outer level of the problem is solved using a derivative-free optimization algorithm, which gives a local optimum, with the starting alignment being the best one produced by a civil engineer.

The article is organized as follows: Subsection 1.1 overviews some of the past research in road design, Section 2 describes basic terminology, Section 3  explains the geometric specifications of a horizontal alignment, Section 4 describes our proposed horizontal alignment optimization models in detail, Section 5 describes the derivative-free optimization solvers we used, Section 6 reports the numerical results for the test problems, and Section 7 summarizes the contributions and highlights some future works.

\subsection{Past Research in Road Design}

As mentioned, road design is commonly divided into three interrelated sub-problems: earthwork optimization, vertical alignment optimization, and horizontal alignment optimization.  Each problem relies on the solution to the sub-problem proceeding it (i.e., vertical alignment optimization requires solving earthwork optimization, and horizontal alignment optimization requires solving vertical alignment optimization).  As such, there is a clear hierarchy in terms of problem difficulty.

Earthwork optimization is perhaps the most established of the three sub-problems.  Many studies investigated earthwork allocation and vertical alignment optimization in road design. Hare et al. \cite{Hare-11} and Lima et al. \cite{LIMA-13} developed two mixed integer linear programming models for earthwork operation in road construction. Unlike Burdett et al. \cite{Burdett-2014} who developed a model for the earthwork allocation problem considering earthwork as discrete 3D blocks, in the present article we use a section-based model, noting that section-based models achieve similar precision as 3D block based model when section lengths are less than 30m~\cite{CHENG-13}.

While vertical alignment optimization is more complicated, it also has a rich research literature.  In 2009 Moreb \cite{MOREB-09} developed a linear programming model combining the vertical alignment and earthwork allocation optimization. In 2010, Koch and Lucet \cite{KOCH-10} advanced Moreb's model by removing unnecessary errors in slope constraints. More recently, Hare et al. \cite{Hare-2014} incorporated the vertical alignment in the earthwork allocation model, resulting in a mixed integer linear programming model that can be solved efficiently in practice.

While horizontal alignment optimization is the most complicated problem, it has nonetheless seen a number of approaches.  Jong et al. \cite{Jong-1998,Jong-2000} developed a horizontal alignment optimization model which was solved by a genetic algorithm. However, the resulting horizontal alignment does not offer any guarantee of (local) optimality.  In 2008, Easa et al. \cite{Easa-2008} developed an optimization model incorporating safety constraints, which were quantified as the expected collisions for an alignment. Although this model guarantees global optimality, the associated vertical alignment cost is not incorporated in the optimization process.  In 2009, Lee et al. \cite{Lee-2009} presented a heuristic based method to optimize the horizontal alignment that works in two stages. In the first stage, the heuristic tries to approximate a piecewise linear alignment and then in the second stage, it refines the solution to make the previously generated piecewise linear alignment compatible with a real road alignment. The solution alignment of the model yields a practical alignment but since a heuristic algorithm was used to solve the model, optimality is not guaranteed \cite{Lee-2009}.

In the literature, a few studies also investigated the problem as a three dimensional alignment optimization problem in which the vertical and horizontal alignments are optimized simultaneously.

Tat and Tao~\cite{Tat-2003} proposed a three-dimensional alignment optimization model, which they solved using a genetic algorithm. Their model considers all of the major constraints in road design. Akay~\cite{Akay-2006} developed a model for three dimensional alignment optimization for forest roads and solved it using a simulated annealing algorithm. Aruga~\cite{Aruga-2005} used a tabu search method to optimize three dimensional alignments of forest roads.

A criteria-based decision support system for three dimensional alignment optimization was developed by Jha \cite{Jha-2003} considering the environmental costs. Jong et al. \cite{Jong-2003} presented an evolutionary model for optimizing the vertical and horizontal alignment simultaneously. The previous two models \cite{Jha-2003,Jong-2003} were improved in \cite{Jha-2006b} by considering accessibility, proximity, and land-use changes, and further improved in \cite{KJSK-2007} to consider incorporating bridge and tunnel costs.

Cheng and Lee \cite{Cheng-2006} also proposed a heuristic-based model for three dimensional alignment optimization. The heuristic solves the models in three steps: first, it generates a good general horizontal alignment by adding, deleting, or moving the intersection points one by one, then it determines an improved horizontal alignment by adjusting the intersection points based on the previously generated horizontal alignment, and finally, it finds a better three dimensional alignment by tuning the vertical alignment corresponding to the previously obtained horizontal alignment.

Most similar to this work, Kang et al. \cite{Kang-2010} developed a bi-level optimization model for road alignment design. In the upper level a set of alternative good alignments is generated and in the lower level an alignment is selected from the alternative alignments obtained in the upper level. The model \cite{Kang-2010} was solved using a genetic algorithm. Recently, Kang et al. \cite{Kang-2012} also proposed a three dimensional alignment optimization model based on genetic algorithm and geographic information system (GIS).

All of the above mentioned three-dimensional alignment optimization models use heuristic-based algorithms which do not guarantee optimality (or even local optimality), and have no or very weak convergence guarantees.

\section{Terminology}
Horizontal alignment optimization consists of finding an optimal curve connecting two given end-points within a designated corridor. The ground profile data in a corridor is given at some discrete points, named \emph{data points}, within the specified corridor.

\begin{figure}[!ht]
\centering
\includegraphics[scale=0.4]{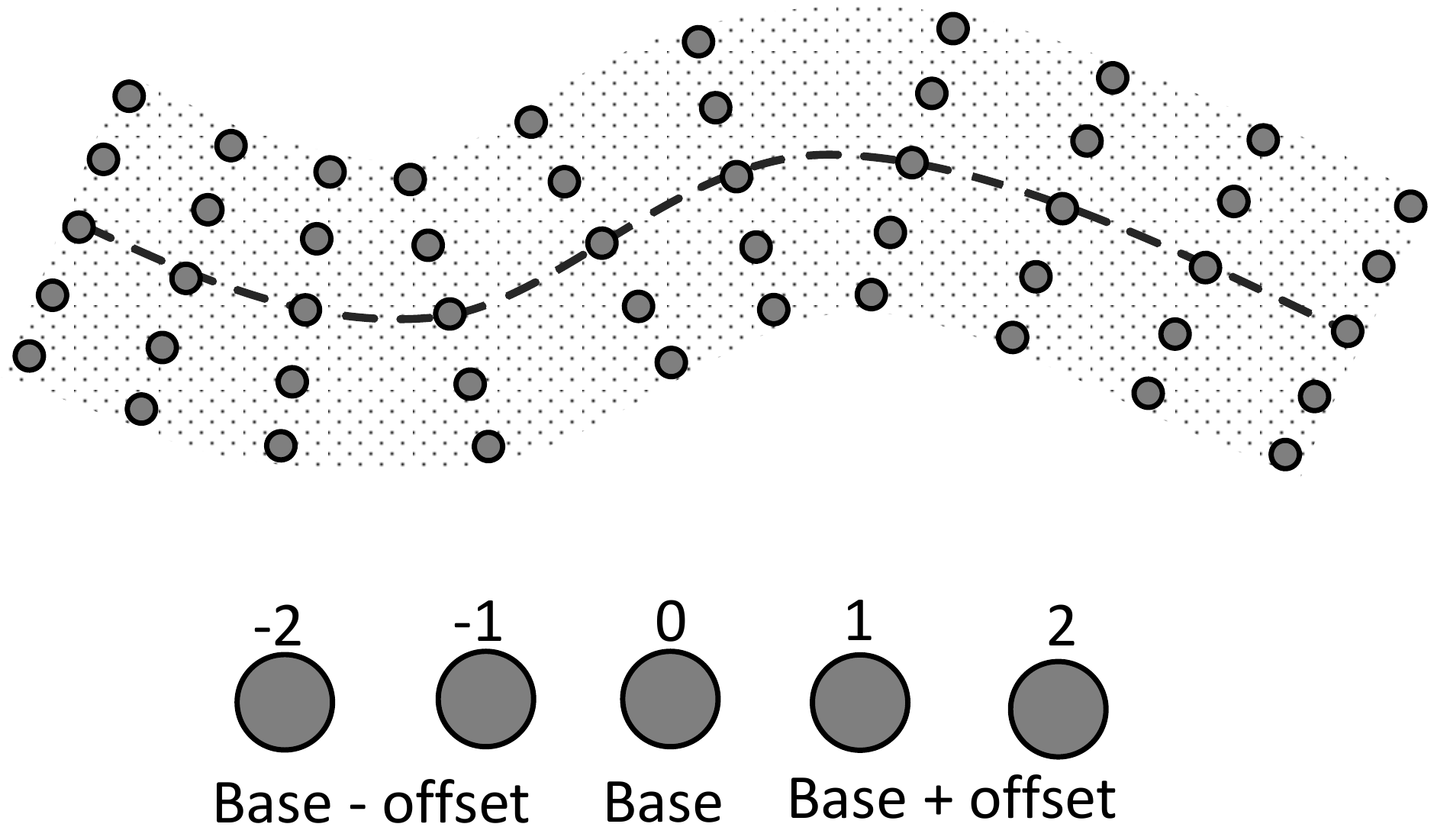}
\caption{Corridor of a horizontal alignment.}
\label{fig:hadata}
\end{figure}

There are two types of data points, namely, \emph{base data points} and \emph{offset data points}. Typically, the base data points are the points along the engineer's original horizontal alignment (however, this is not strictly necessary). The offset data points represent the horizontal displacement from the base data points.
The base data points are selected a few units apart between the two end-points along the baseline. Each of the base data points has some associated offset data points in both the left and the right directions, see Figure~\ref{fig:hadata}. A base data point together with the associated offset data points is defined as a \emph{station}. The baseline of a corridor is a curve connecting the base data points (i.e., the dotted curve in Figure~\ref{fig:hadata}). In practice, a baseline is a primarily defined alignment by engineers.

Base data points, offset data points, stations, and baselines are all fixed input data for a given horizontal alignment problem.  Ultimately, the horizontal alignment problem is to determine the optimal curve (in terms of vertical alignment and earthwork cost) through those stations, subject to road design constraints.

\begin{figure}[!ht]
\centering
\includegraphics[scale=0.4]{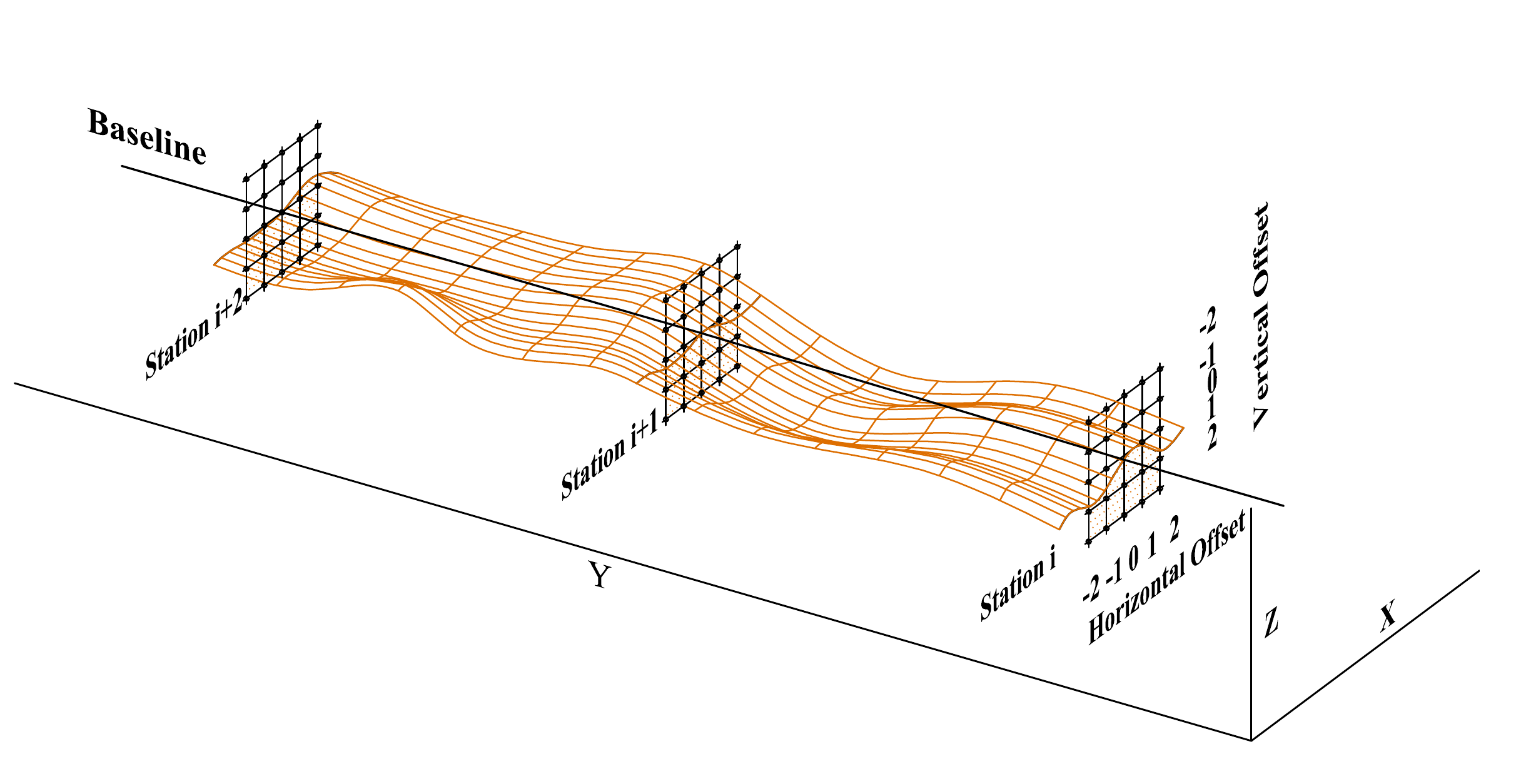}
\caption{Three dimensional corridor. Corresponding to each horizontal offset, a set of vertical offsets is shown along the vertical line}
\label{fig:3dalignment}
\end{figure}

Each data point within the corridor, either a base data point or an offset data point, has some associated ground profile data. Therefore, for the vertical road profile, we can move vertically up and down for each horizontal data point. The horizontal and the vertical displacements from the baseline make a discrete grid for each station, see Figure~\ref{fig:3dalignment}. Our goal is to find, for each station, a horizontal offset that generates a horizontal alignment and a vertical alignment which is (locally) optimal.

\section{Geometric representation}
\begin{figure}[!ht]
\centering
\includegraphics[scale=0.3]{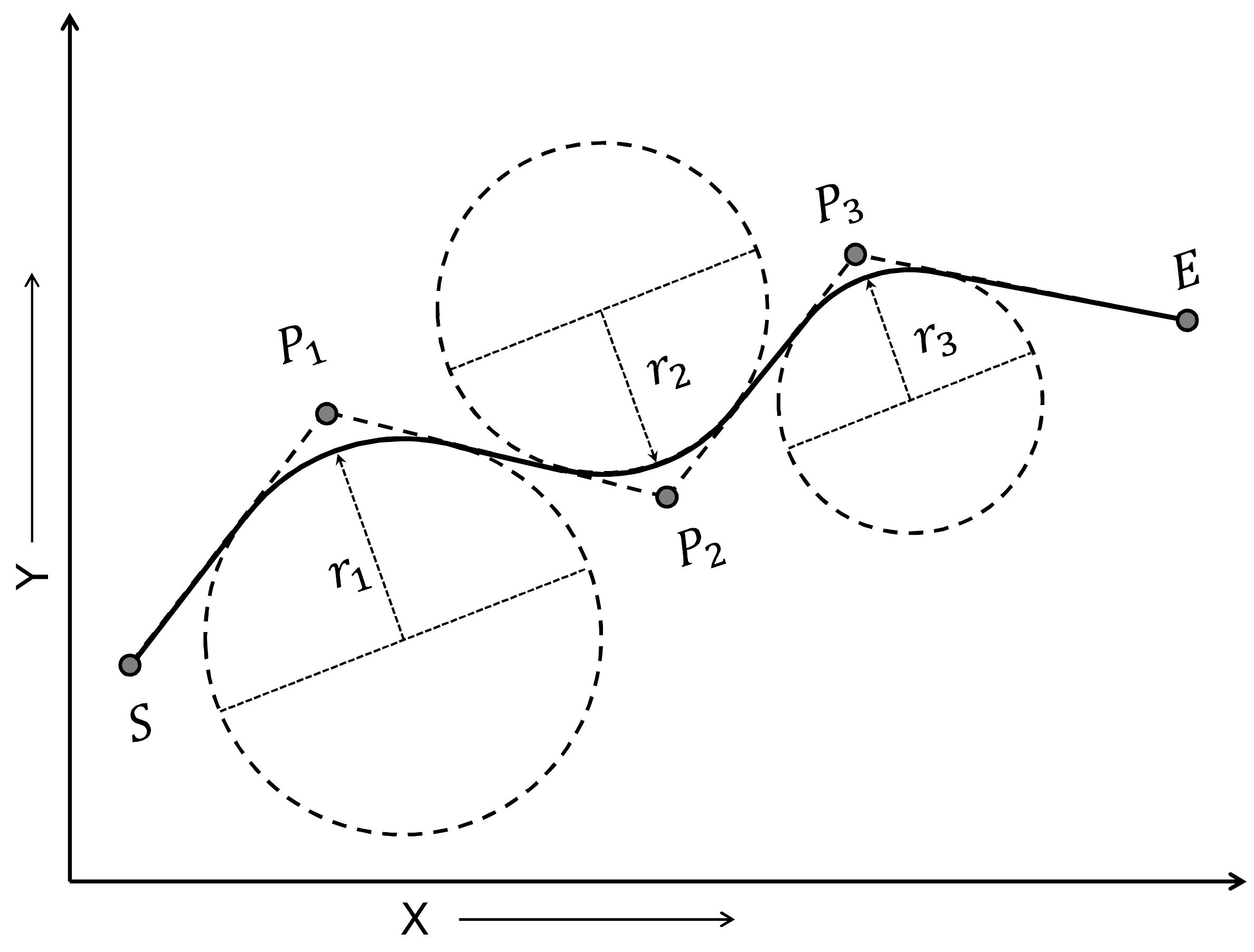}
\caption{Geometric representation of an horizontal alignment.}
\label{fig:geometric_representation_HA}
\end{figure}

A horizontal alignment consists of a sequence of circular curves and tangential lines, which are defined by some intersection points and the radius of curvature associated with each intersection point. In Figure~\ref{fig:geometric_representation_HA}, $S$ and $E$ are the start and end points of the alignment, respectively. The intersection points of the alignment are $P_1$, $P_2$, and $P_3$. Each intersection point has a radius of curvature that defines the circular curve. The radius of curvature associated with the intersection points $P_1$, $P_2$, and $P_3$ are $r_1$, $r_2$, and $r_3$, see Figure~\ref{fig:geometric_representation_HA}.

Let $i$ be the index of the intersection points and $n$ be the number of intersection points. Since the intersection point $P_i$ has an associated radius of curvature $r_i$, we define an intersection point with radius of curvature as $(P_i,r_i)$, where $P_i \in \mathbb{R}^2$  and  $r_i \in \mathbb{R}$.
Without loss of generality, we can say that the start and end points are two points in $\mathbb{R}^2$ with zero radius of curvature and note them $(P_0,0)$ and $(P_{n-1},0)$. So we represent a horizontal alignment $HA$ as the sequence
\begin{equation}
HA=((P_0,0),(P_1,r_1),(P_2,r_2),\hdots,(P_{n-1},0)).
\end{equation}

\begin{figure}[!ht]
\centering
\includegraphics[scale=0.3]{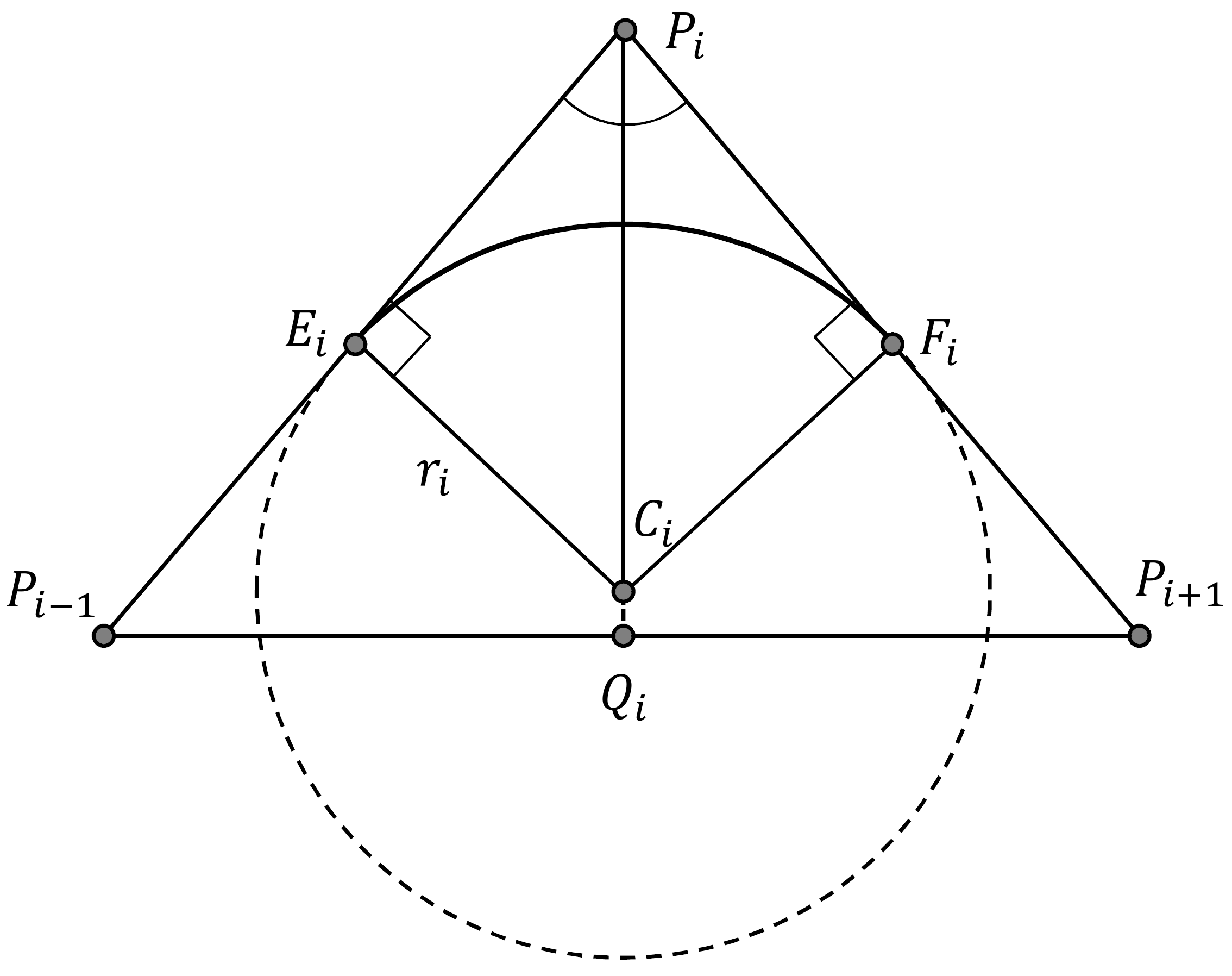}
\caption{Geometric specifications of a circular curve.}
\label{fig:gemetric_specification_circular_curve}
\end{figure}

To determine the actual horizontal alignment, we need to calculate the circular curves and tangential line segments from the given intersection points  and the associated radius of curvatures.
Let $E_i$ and $F_i$ be the left and right tangential point, respectively, see Figure~\ref{fig:gemetric_specification_circular_curve}.
Let $C_i$ be the center of curvature that corresponds to the intersection point $P_i$.

\subsection{Tangential point calculation}
We define the following variables: $\theta_i$ is the angle at $P_i$ using three consecutive intersection point $P_{i-1}$, $P_i$, and $P_{i+1}$, $Q_i$ is the intersection point of the angle bisector of $\theta_i$ and the line joining $P_{i-1}$ and $P_{i+1}$, $\textbf{U}_i =P_{i-1}-P_i$, $\textbf{V}_i =P_{i+1}-P_i$, and $\textbf{W}_i =P_{i+1}-P_{i-1}$. We have
\begin{equation}
\theta_i= \arccos \left(\frac{\textbf{U}_i\cdot\textbf{V}_i}{\norm{\textbf{U}_i} \norm{\textbf{V}_i}} \right).
\label{eq:theta_calcluation}
\end{equation}
Since the angle bisector $P_iQ_i$ bisects the angle $\theta_i$, we have $\angle P_{i-1}P_iQ_i=\frac{\theta_i}{2}=\angle P_{i+1}P_iQ_i$. The lines $P_iP_{i-1}$ and $P_iP_{i+1}$ are the tangent to the circle at the tangential point $E_i$ and $F_i$. Thus we have $C_iE_i\bot P_iP_{i-1}$ and $C_iF_i\bot P_iP_{i+1}$. The triangle $\triangle P_iE_iC_i$ and $\triangle P_iF_iC_i$ are right-angle triangles. The segment $P_iC_i$ is the common side of $\triangle P_iE_iC_i$ and $\triangle P_iF_iC_i$. Since $E_iC_i$=$F_iC_i$, we have $P_iE_i=P_iF_i$. Let $l_t$ be the length of $P_iE_i$. We have
\begin{equation}
l_t = \frac{r_i}{\tan \frac{\theta_i}{2}}.
\end{equation}
Let $\hat{e}_{\textbf{U}_i}$ and $\hat{e}_{\textbf{V}_i}$ be the two unit vector of $\textbf{U}_i$ and $\textbf{V}_i$.  The tangential point $E_i$ and $F_i$ is calculated as
\begin{align}
E_i&=P_i+l_t\hat{e}_{\textbf{U}_i} \label{eq:left_tangentail_point},\\
F_i&=P_i+l_t\hat{e}_{\textbf{V}_i} \label{eq:right_tangentail_point}.
\end{align}

\subsection{Center point calculation}

\begin{fact}[{\cite[Book \rom{6} Proposition \rom{3}]{Byrne-87}}]
\label{angel_bisector_theorem}
The angle bisector of an angle in a triangle divides the opposite side in the same ratio as the sides adjacent to the angle.
\end{fact}

Let $l_b$ be the length of $Q_iP_{i-1}$. The length of $P_iP_{i-1}$, $P_iP_{i+1}$ and $P_{i-1}P_{i+1}$ are $\norm{\textbf{U}_i}$, $\norm{\textbf{V}_i}$, and $\norm{\textbf{W}_i}$, respectively.
The line $P_iQ_i$ is the angle bisector of  $\theta_i$ in $\triangle P_{i-1}P_iP_{i+1}$. Using Fact~\ref{angel_bisector_theorem} we have
\begin{equation}
l_b = \frac{\norm{\textbf{U}_i} \norm{\textbf{W}_i}}{\norm{\textbf{U}_i}+\norm{\textbf{V}_i}}.
\end{equation}
Let $\hat{e}_{\textbf{W}_i}$ be the unit vector of  $\textbf{W}_i$. The point $Q_i$ can be calculated as
\begin{equation}
Q_i= P_{i-1}+ l_b \hat{e}_{\textbf{W}_i}.
\end{equation}
We note $\textbf{X}_i= Q_i-P_i$. Let $l_x$ be the length of $P_iC_i$. We have
\begin{equation}
l_x =\frac{l_t}{\cos \frac{\theta_i}{2}}.
\end{equation}
Let $\hat{e}_{\textbf{X}_i}$ be the unit vector of  $\textbf{X}_i$.
So the center point is obtained by
\begin{align}
C_i&=P_i+l_x\hat{e}_{\textbf{X}_i}. \label{eq:center_point}
\end{align}

\section{Model description}
We formulate the horizontal alignment optimization problem as a bi-level optimization problem. The inner problem is the vertical alignment optimization problem, while the outer problem uses this information to seek an optimal horizontal alignment.

\subsection{Variables and objective function}
\begin{figure}[!ht]
\centering
\includegraphics[scale=0.4]{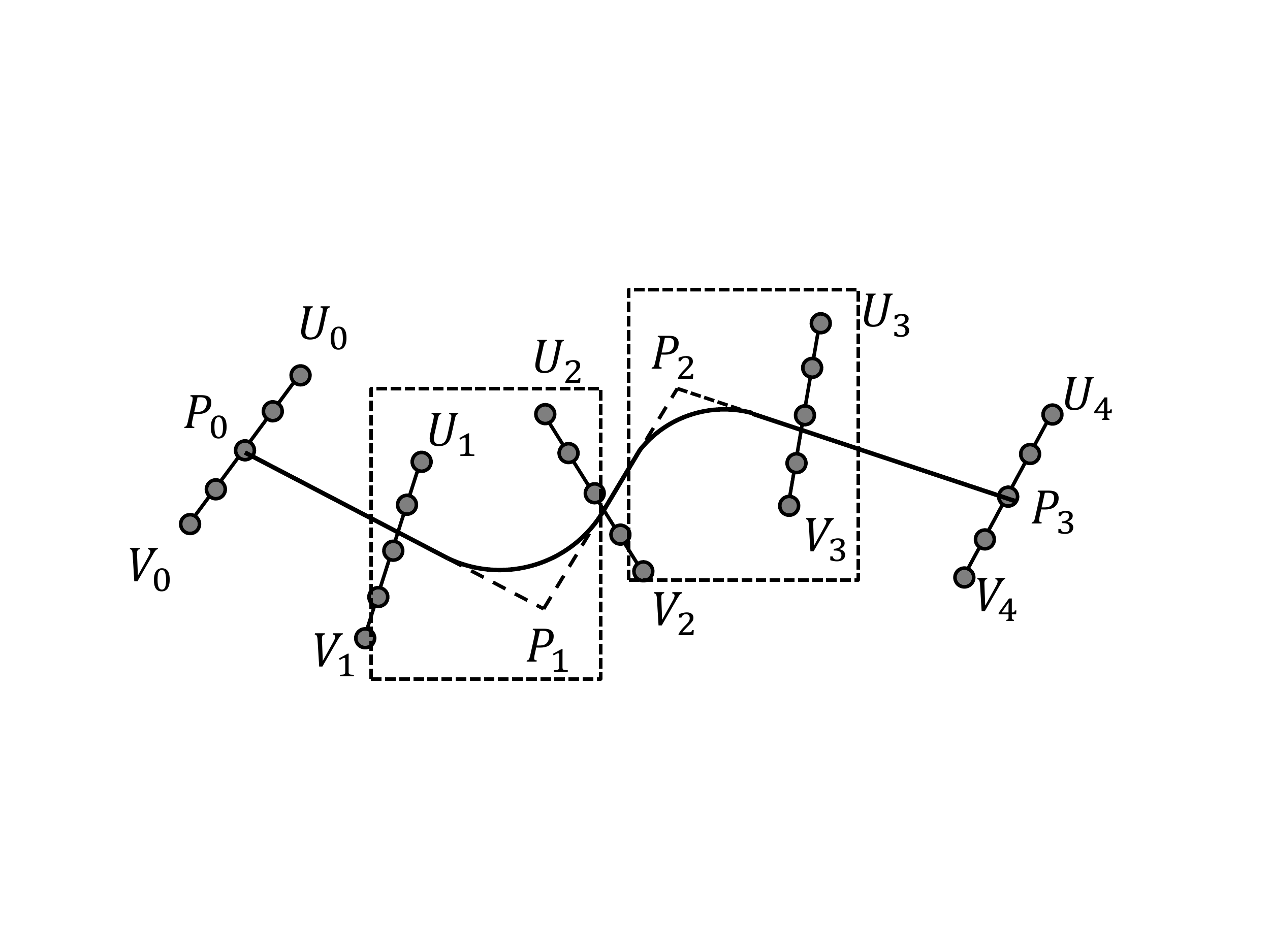}
\caption{A potential horizontal alignment in a specified corridor.}
\label{fig:HA_in_a_specified_corridor}
\end{figure}

We can define a horizontal alignment in a specified corridor using a set of intersection points and the associated radius of curvatures, see Figure~\ref{fig:HA_in_a_specified_corridor}. By varying the location of the intersection points in the associated feasible regions, a wide variety of horizontal alignments can be built. Note that in our model the feasible region of an intersection point is defined as a rectangular box.  Since the intersection points of a horizontal alignment are in the \emph{xy}-plane, we denote $P_i=(x_i,y_i)$. We write the variable vector as
\begin{equation}\begin{array}{rl}
X&={((x_1,y_1,r_1),(x_2,y_2,r_2)\hdots,(x_{n-2},y_{n-2},r_{n-2}))}\\
&={((P_1,r_1),(P_2,r_2)\hdots,(P_{n-2},r_{n-2}))}.
\end{array}
\end{equation}
Note that the starting point $P_0$ and the ending point $P_{n-1}$ are fixed and do not have any corresponding radius variables.

Suppose we have $n_s$ stations indexed from 0 to $n_s-1$. At each station, we define a cross-section line which passes through the leftmost offset data point and the rightmost offset data point. Let $U_j$ and $V_j$ be the leftmost and rightmost offset data points of a station (see Figure \ref{fig:HA_in_a_specified_corridor}). So the cross-section lines in a specified corridor are written as
\begin{equation}
\label{eq:cross_section_line}
\begin{aligned}
	L_j(t) = (1-t)U_j+tV_j \;\;\;\;\;\;\;\; \text{for } t \in \mathbb{R},\mbox{ } j=0,1,\hdots ,n_s-1.
\end{aligned}
\end{equation}

Since a horizontal alignment is a linear-circular curve, we can calculate the point of intersections between the cross-section lines and a horizontal alignment. A point of intersection can be represented as a parameter value of the cross-section line which is denoted as $t_j$. If $t_j \in [0,1]$ for all $j$, then a horizontal alignment is inside of the corridor. The values of $t_j$ can be mapped into the horizontal offsets. In our input data, horizontal offsets are discrete values, see Figure~\ref{fig:3dalignment}. So we used linear interpolation to compute vertical ground profile data for any (continuous) offset values. The accuracy of this interpolation is clearly dependent on the distance between the offset data points. If higher accuracy is desired, then more offset data points should be included in the problem formulation.  As computation time of this interpolation is independent of the number of offset data points, increasing the number of data points should only effect the data storage requirements, not solution time.

Once we have vertical ground profile data at each station for a horizontal alignment we can solve a vertical alignment optimization problem using the models developed in \cite{Hare-2014}.  Since we will use the vertical alignment solution, denoted $\mathcal{C_{VA}}$, as a black-box that we seek to minimize, it is not necessary to understand how exactly $\mathcal{C_{VA}}$ is computed for this paper.  Nonetheless, in Appendix \ref{appendix}, we provide a brief overview of the mixed-integer linear program used to compute $\mathcal{C_{VA}}$ (for curious readers).

The final optimization problem can be written as
\begin{equation}\label{problem}
{\min f(X)=}
\left\{
\begin{array}{ll}
\mathcal{C_{VA}}(X)  & \mbox{if } t_j\in \left[0,1\right], \mbox{ } j=0,1,\hdots ,n_s-1. \\
\infty & \mbox{ otherwise, }
\end{array}\right.
\end{equation}
where $\mathcal{C_{VA}}(X)$ is a function that gives the optimal vertical alignment cost for $X$.   Note that when a horizontal alignment is outside of the corridor, we set the cost to infinity.

\subsection{Constraints}
\begin{figure}[!ht]
\centering
\includegraphics[scale=0.2]{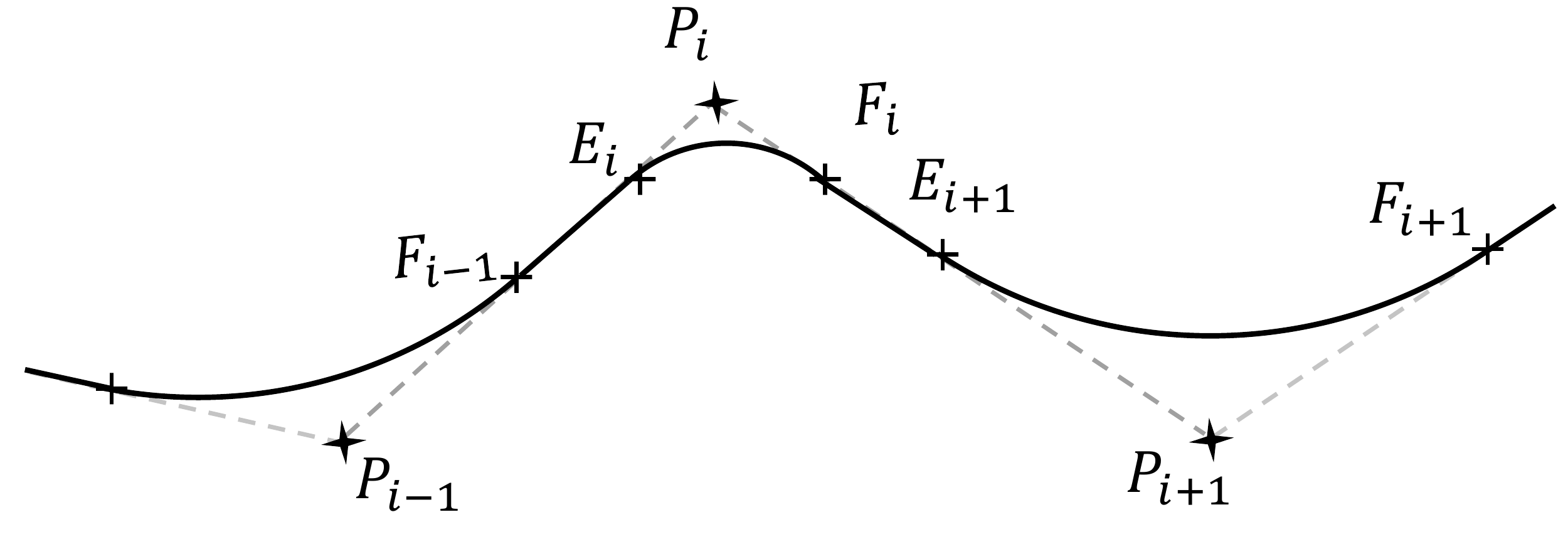}
\includegraphics[scale=0.2]{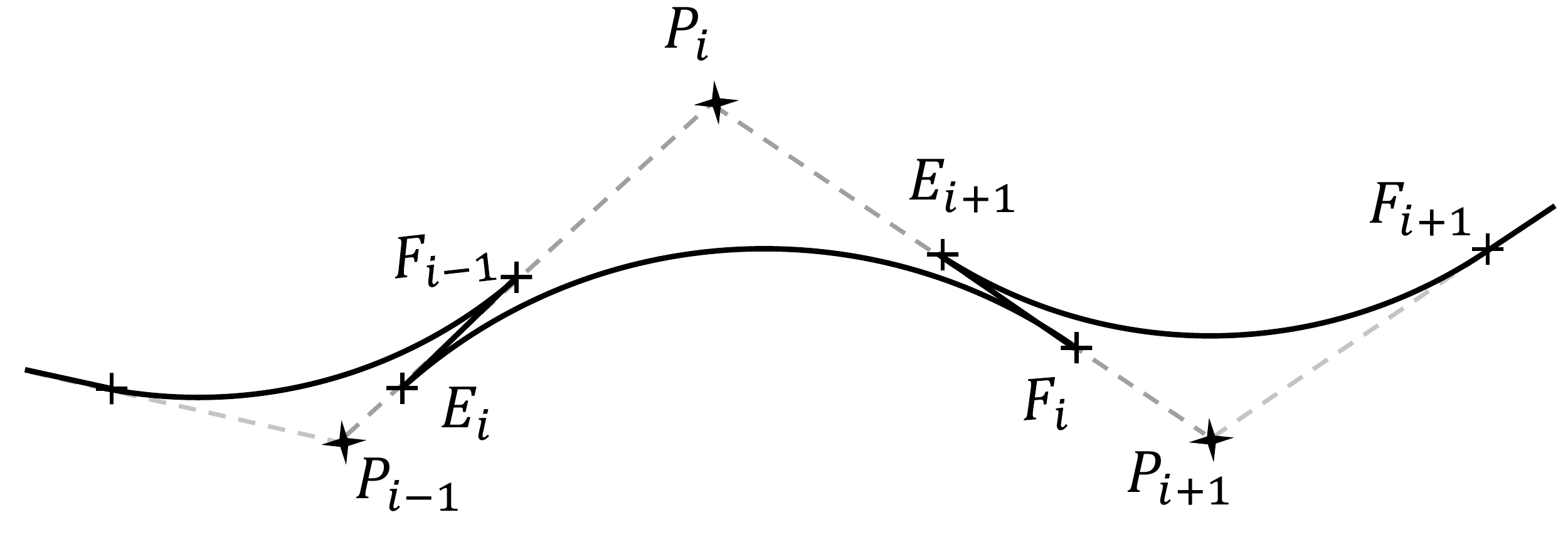}
\caption{Left: a curve with well-defined tangent lines.  Right: a curve without well-defined tangent lines.}
\label{fig:discontinuty_in_HA_segment}
\end{figure}

A horizontal curve consists of tangential line segments followed by circular arcs. Two consecutive circular arcs are connected by a tangential line.  In order for this tangent line to be well-defined, it is important that the circular arcs are small enough that they do not `overlap'.  In Figure~\ref{fig:discontinuty_in_HA_segment}, we give an pictorial example of what goes wrong if the circular arcs are too large.
Therefore, in order to maintain well-defined tangential line segments, on the line passing through the intersection points $P_{i-1}$ and $P_i$, the length of $P_{i-1}P_{i}$ must be greater than or equal to the summation of the length of $P_{i-1}F_{i-1}$ and the length of $P_{i}E_{i}$.
We can write these constraints as
\begin{equation}
\norm{P_{i}-P_{i-1}} \geq \norm{P_{i-1}-F_{i-1}}+\norm{P_{i}-E_{i}}, \;\;\; i=1,2,\hdots,n-2.
\end{equation}
We will refer to these as the continuity constraints.

If the radius of curvature is too small (or zero) then a horizontal alignment might get a sharp turn. So the optimal radius of curvature must be greater than or equal $R_{\min}$, the minimum radius of curvature. In other words,
\begin{equation}
r_i\geq R_{\min}, \;\;\; i=1,2,\hdots,n-2.
\end{equation}

In our model, each intersection point has a feasible region defined by a rectangular box, see Figure~\ref{fig:HA_in_a_specified_corridor}. A given rectangular box is defined using the leftmost bottom corner point $({l_x}_i, {l_y}_i)$ and the rightmost top corner point $({u_x}_i, {u_y}_i)$. Thus the box constraints corresponding to intersection point $P_i$ are written as

\begin{equation}
{l_x}_i\leq x_i \leq {u_x}_i \mbox{ and } {l_y}_i\leq y_i \leq {u_y}_i,\mbox{  } i= 1,2,\hdots,n-2.
\end{equation}

\section{Solution Approach}
The objective function of problem \eqref{problem} is an optimization problem in itself.  Indeed, $\mathcal{C_{VA}}$ is evaluated through a large scale mixed integer linear program (see Appendix \ref{appendix}). As such, it is very hard to access the derivative information of the objective function (if it exists).  As a result, optimizing the objective function cannot be accomplished by gradient based methods (such as BFGS) or structure based methods (such as the Simplex algorithm).  Instead, it is necessary to apply some non-gradient based approach.  We specifically turn our attention to two derivative-free optimization (DFO) solvers: NOMAD \cite{NOMADMan-2009} (version 3.5, available at \url{http://www.gerad.ca/nomad}) and HOPSPACK \cite{Hops-2009} (version 2.0.2, available at \url{http://www.sandia.gov/hopspack}).  In Subsections \ref{ss:NOMAD} and \ref{ss:HOPSPACK} we give very brief overviews on how each of these solvers work.  For greater detail, we refer readers to \cite{NOMADMan-2009} (for NOMAD) and \cite{Hops-2009} (for HOPSPACK).

\subsection{NOMAD}\label{ss:NOMAD}

NOMAD abbreviates ``Nonlinear Optimization with MADS''.  The solver is an freely available and based on the ``Mesh Adaptive Direct Search'' (MADS) algorithm designed in \cite{Audet-2006}.   MADS falls into the large class of {\em pattern search} methods, which essentially iterate the following.
    \begin{quote} Given an iterate $x^k$, a search radius $\Delta^k>0$, and a `pattern' $\{d^1, d^2, ... d^m\}$ with $\max \{ |d^i| : i=1, 2, ...,m \} = \Delta^k$, evaluate $f(x^k + d^i)$ for each $i=1, 2, ... m$.  If a better point $x^k + d^j$ is found, then set $x^{k+1} = x^k + d^j$ and repeat using $\Delta^{k+1}=\Delta^k$.  If no better point is found, then set $\Delta^{k+1}=\gamma\Delta^k$ where $0 < \gamma < 1$ and repeat.
    \end{quote}
Pattern search methods enjoy several strong feature.  First, pattern searches only require function evaluations, not gradient evaluations, to optimize the objective function.  Second, unlike more heuristic methods (like genetic algorithms) in many cases pattern searches can be mathematically proven to converge to a local minimizer. Third, pattern searches are naturally parallelized.

The MADS algorithm follows specific mathematical rules on how to generate the pattern, and using these rules can guarantee convergence to a local minimizer \cite{Audet-2006}.  While the rules are complex, they are nicely described (with examples) in \cite{Audet-2006}. 

NOMAD is deterministic in nature, meaning that if the algorithm is run twice on the same problem, it will return the exact same outcome both times.  As such, in our numerical tests, NOMAD is only run once per test problem.

The NOMAD solver allows for several different stopping criterion.  For our numerical tests, we employ the `minimum mesh size' stopping criterion, which (in the language above) stops the algorithm when the search radius ($\Delta^k$) becomes too small.

\subsection{HOPSPACK}\label{ss:HOPSPACK}

HOPSPACK, the Hybrid Optimization Parallel Search PACKage, implements the Asynchronous Parallel Pattern Search (APPS) algorithm in a C++ software framework \cite{Kolda-2005,Gray-2006}.  Similar to MADS, the APPS algorithm is a pattern search method.  Although the APPS algorithm uses a different set of rules for pattern generation, it is also proven to converge to a local minimizer \cite{Kolda-2005}.  However, unlike MADS and APPS, the HOPSPACK algorithm includes some heuristic techniques to break free of local minimizers.  These heuristics invalidate the original convergence proof of APPS, but have been found effective in practice.  (If the heuristic are only employed a finite number of times, then the convergence proof of APPS still holds, but it is unclear if this safe-guard is implemented in HOPSPACK.)  The heuristics are also stochastic in nature, so HOPSPACK will not necessarily return the same optimized outcome if run a second time on the same problem.  Therefore, in our numerical tests, we have run HOPSPACK five times for each test problem.

Similar to the MADS algorithm, HOPSPACK includes a stopping criterion based on a minimum allowable search radius (this time called `minimum step size'). For consistency with NOMAD, in our numerical tests we use this stopping criterion.

\section{Experimental results}
We performed numerical experiments on five different roads listed in Table~\ref{table:test_problem_set}. The road profile data correspond to real roads and were given by our industry partner Softree Technical System Inc.

\begin{table}[!ht]
\centering
\caption{Specifications of the test problems}
\begin{tabular}{|p{2cm} | >{\centering\arraybackslash} p{3 cm}| >{\centering\arraybackslash} p{3 cm}| }
\hline
Road Name & Length of the Road (meters) & Number of Intersection Points \\
\hline \hline
Road A & 358 & 8  \\
Road B & 4000 & 5 \\
Road C & 5013 & 14 \\
Road D & 1156 &  22 \\
Road E & 721 &  9 \\
\hline
\end{tabular}
\label{table:test_problem_set}
\end{table}
All of the experiments were carried out in a Dell workstation with an Intel(R) Xenon(R) 2.40 GHz (2 cores) processor, 24 GB of RAM and a 64-bit Windows 7 Enterprise operating system. Software to process and solve the optimization problem was implemented in C++ using Microsoft Visual Studio 2010 Professional Edition.

Since we are interested in a solution close to a local minimum, for both solvers (NOMAD and HOPSPACK) we set the stopping criterions (i.e., minimum mesh size and minimum step length) to 0.1.  As the input data of our model are given in meters, the final scaling of the variables of our model goes down below 10 cm, which means a local minimum should exist in less than 10 cm distance (see \cite{NOMADMan-2009} and \cite{Hops-2009}).  In all tests, convergence to this level of local minimizer was confirmed by the algorithm.

We used the baseline alignment of a corridor as an initial starting point for both solvers. Recall, the NOMAD solver gives a deterministic solution (i.e., different independent runs of the algorithm yield the same solution), while the HOPSPACK solver gives a non-deterministic solution (i.e., different independent runs of the algorithm might yield different solutions). So first, we solved the test problems using the NOMAD solver and then compare with the HOPSPACK solver solutions obtained by several independent runs.

\subsection{Case Study}

\begin{figure}[!ht]
\centering
\includegraphics[scale=0.35]{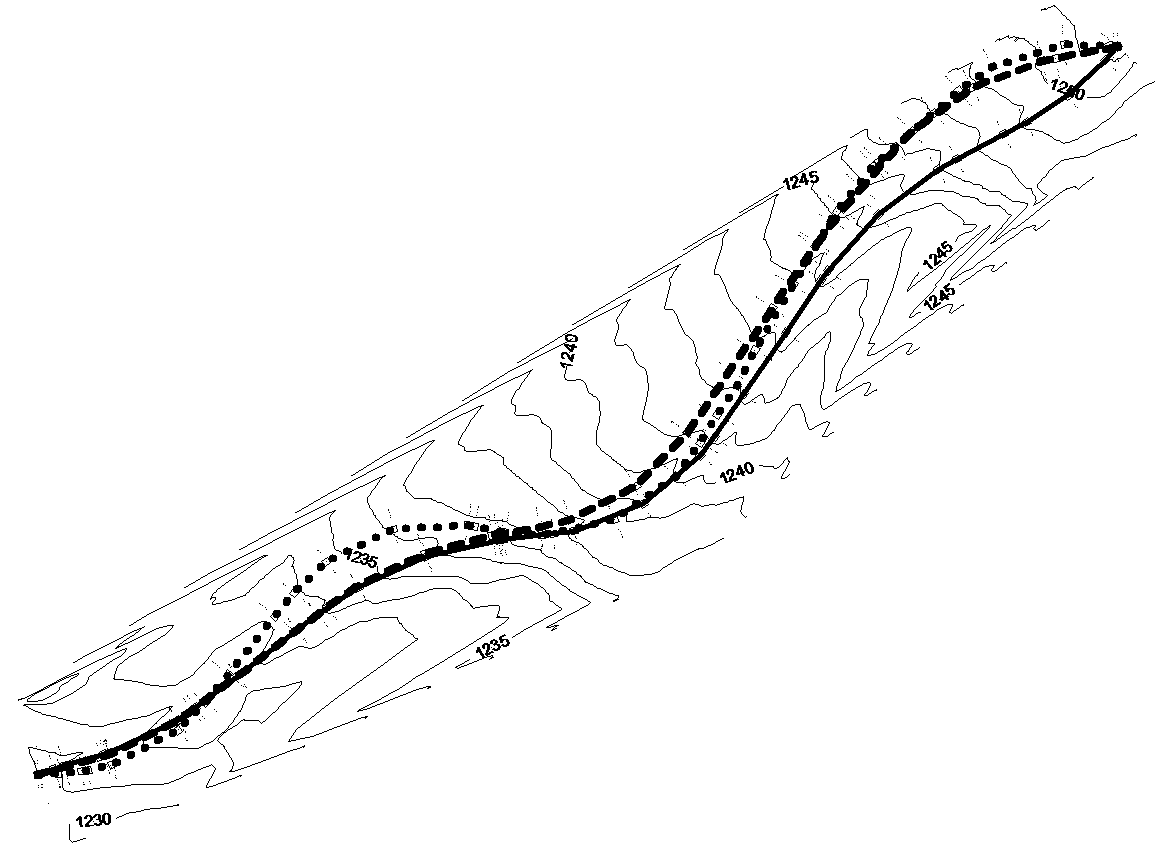}
\caption{Initial and Improved Horizontal Alignments for Road A. \newline
    {\small \tt Solid Black = Initial Alignment, Dashed = NOMAD optimized alignment, Dotted = HOPSPACK optimized alignment}}
\label{fig:roadA}
\end{figure}

Before providing a summary of the results concerning all test problems, we focus on a single problem as a case study.  Specifically, we discuss Road A, which is a short (358m) section of a forestry road.  In Figure \ref{fig:roadA}, we visualize the initial baseline alignment, and the two optimized alignments (as solved by NOMAD and HOPSPACK).  For the sake of visualization, Figure  \ref{fig:roadA} also includes the contour lines of the terrain within the allowable corridor for Road A.

In Figure \ref{fig:roadA}, the initial baseline alignment for Road A is presented as a solid black line.  Notice it goes exactly down the middle of the corridor.  (This is not a required feature of the initial alignment, it is simply because the test problem data, provided by Softree Inc., happened to take this form.)  The optimal vertical alignment for the baseline alignment of Road A cost \$1897.

In Figure \ref{fig:roadA}, the optimized horizontal alignment for Road A, when computed using NOMAD as the DFO solver is presented as a dashed line.   Recall that HOSPACK in non-deterministic, and so HOPSPACK was run five times on each test problem.  In Figure \ref{fig:roadA}, the first of these five optimized horizontal alignments when computed using HOPSPACK is presented as a dotted line.  The optimal vertical alignment for the horizontal alignment created using NOMAD costs \$1361, and the optimal vertical alignment for the horizontal alignment created using HOPSPACK costs \$1291.  Both of these represent significant savings over the original cost.  In particular, the percentage savings can be computed as
    \[\frac{1897-1361}{1897}100\% = 28.3\% \quad \mbox{and} \quad \frac{1897-1291}{1897}100\% = 31.9\%.\]

\subsection{Summary results for all test problems}
We now turn our attention to an aggregate analysis of the all test problems.

Table~\ref{table:test_problems_cost_imporvement_NOMAD} shows the cost improvement of the objective functions, the number of black-box evaluation and wall-clock time required to solve the test problems using the NOMAD solver. Note that a candidate alignment that falls outside the corridor counts as a function call but does not require solving a MILP, which explains that road C and D require similar time to solve with Road~D having double the number of function calls. The cost improvements range from 11$\%$ to 54$\%$ with an average of $27\%$.  While the wall-clock times for some problems are high, the construction of a new road can take years.  As such, the used computer time (which could be run over night) is not unreasonable given the average saving in final costs.
\begin{table}[!ht]
\caption{Cost improvement, no. of black-box evaluations and wall-clock time required to solve the test problems by the NOMAD solver.}
\centering
\begin{tabular}{|L{1.3cm}| R{1.7 cm}| R{1.7 cm}| R{1.5 cm}| R{1.5 cm}| R{1.5 cm}|}
\hline
Road Name & \multicolumn{1}{>{\centering\arraybackslash}m{1.7cm}|}{Initial alignment cost}  & \multicolumn{1}{>{\centering\arraybackslash}m{1.7cm}|}{Optimized alignment cost} & \multicolumn{1}{>{\centering\arraybackslash}m{1.5cm}|}{Cost Improv- ement ($\%$)} & \multicolumn{1}{>{\centering\arraybackslash}m{1.5cm}|}{No. of Black-box evaluations} & \multicolumn{1}{>{\centering\arraybackslash}m{1.5cm}|}{Wall-clock time (seconds)} \\
\hline \hline
 Road A & 1,897 & 1,361 & $28\%$ & 2,073 & 1,445 \\
 Road B & 17,036 & 15,198 & $11$\% & 2,528 & 1,770 \\
 Road C & 87,829 & 69,621 & $21 \%$ & 37,165 & 45,647 \\
 Road D & 31,031 &  14,418 & $54\%$ & 90,535 & 47,613 \\
 Road E & 8,054 &  6,498 & $19 \%$ & 11,101 & 5,588 \\
 \hline
\end{tabular}
\label{table:test_problems_cost_imporvement_NOMAD}
\end{table}

Table~\ref{table:comparison_objecitve_functions_values_test_problems} lists the optimum values of the objective functions and the number of black-box evaluations required to obtain the solution for five independent executions of the HOPSPACK solver for each test problem. The differences in the optimum objective function values are calculated with respect to the value obtained by the NOMAD solver. So in  Table~\ref{table:comparison_objecitve_functions_values_test_problems}, a $``+"$ value in the \emph{Diff. in costs} column (i.e., difference in optimum costs obtained by the two solvers) indicates the HOPSPACK solver yields a better solution than the NOMAD solver and a $``-"$ value indicates the opposite.
Similarly, we also calculated the difference in number of black-box evaluations with respect to the number of black-box evaluations required by the NOMAD solver. So in Table~\ref{table:comparison_objecitve_functions_values_test_problems}, a $``+"$ value the \emph{Diff. in  black-box eval} column (i.e., differences in no. of black-box evaluations) indicates the HOPSPACK solver required less black-box evaluations than the NOMAD solver and a $``-"$ value indicates the opposite.

\begin{table}[!ht]
\centering
\caption{Comparison of optimum objective function values  and required no. of black-box evaluations for the HOPSPACK and NOMAD solvers to solve the test problems.}
\begin{tabular}{|L{.5 cm} | C{.5 cm} |  R{1.2 cm}|  R{1.2 cm}|  R{1.1 cm}| R{1.2 cm}|  R{1.2 cm}|  R{1.2 cm}|}
\hline
Rd & Test run & \multicolumn{1}{>{\centering\arraybackslash}m{1.2 cm}|}{Opt. cost-NOMAD}   & \multicolumn{1}{>{\centering\arraybackslash}m{1.2 cm}|} {Opt. cost-HOPS.} & \multicolumn{1}{>{\centering\arraybackslash}m{1 cm}|} {Diff. in costs (\%)}  & \multicolumn{1}{>{\centering\arraybackslash}m{1.2 cm}|}{$\#$ black-box eval-NOMAD}   & \multicolumn{1}{>{\centering\arraybackslash}m{1.2 cm}|} {$\#$ black-box eval-HOPS.} & \multicolumn{1}{>{\centering\arraybackslash}m{1.2 cm}|} {Diff. in  black-box eval ($\%$)} \\
\hline \hline
& $\#$1 & 1,361 & 1,291  & +5.2$\%$ & 2,073 & 325  & +84.3$\%$\\
& $\#$2 & 1,361 & 1,423  & -4.6$\%$ & 2,073 & 336 & +83.8$\%$\\
A & $\#$3 & 1,361 & 1,418  & -4.2$\%$ & 2,073 & 697 & +66.4$\%$\\
& $\#$4 & 1,361 & 1,278  & +6.1$\%$ & 2,073 &  486 & +76.6$\%$\\
& $\#$5 & 1,361 & 1,486  & -9.2$\%$ & 2,073 &  665 & +67.9$\%$\\
\hline
& $\#$1 & 15,198 & 15,510  & -2.1$\%$ & 2,528 & 316  & +87.5$\%$\\
& $\#$2 & 15,198 & 15,141  & +0.4$\%$ & 2,528 & 309 & +87.8$\%$\\
B & $\#$3 & 15,198 & 15,128  & +0.5$\%$ & 2,528 & 286 & +88.8$\%$\\
& $\#$4 & 15,198 & 15,172  & +0.2$\%$ & 2,528 &  392 & +88.5$\%$\\
& $\#$5 & 15,198 & 15,529  & -2.2$\%$ & 2,528 &  485 & +80.8$\%$\\
\hline
& $\#$1 & 69,621 & 70,161  & -0.8$\%$ & 37,165 & 7,547  & +79.7$\%$\\
& $\#$2 & 69,621 & 70,378  & -1.1$\%$ & 37,165 & 31,418 & +15.5$\%$\\
C & $\#$3 & 69,621 & 69,995  & -0.5$\%$ & 37,165 & 3,392 & +90.9$\%$\\
& $\#$4 & 69,621 & 67,301  & +3.3$\%$ & 37,165 &  3,213 & +91.4$\%$\\
& $\#$5 & 69,621 & 67,045  & +3.7$\%$ & 37,165 &  4,661 & +87.5$\%$\\
\hline
& $\#$1 & 14,418 & 13,190 & +8.5$\%$ & 90,535 & 15,852  & +82.5$\%$\\
& $\#$2 & 14,418 & 15,154  & -7.6$\%$ & 90,535 & 19,997 & +77.9$\%$\\
D & $\#$3 & 14,418 & 14,155 & +1.8$\%$ & 90,535 & 21,194 & +76.6$\%$\\
& $\#$4 & 14,418 & 14,016  & +2.8$\%$& 90,535 & 17,816 & +80.3$\%$\\
& $\#$5 & 14,418 & 15,384  & -6.7$\%$& 90,535 & 21,339 & +76.4$\%$\\
\hline
& $\#$1 & 6,497 & 6,524  & -0.4$\%$ & 11,101 & 2,019  & +81.8$\%$\\
& $\#$2 & 6,497 & 6,475  & +0.3$\%$ & 11,101 & 4,332 & +60.1$\%$\\
E & $\#$3 & 6,497 & 6,497 & 0.0$\%$ & 11,101 & 1,996 & +82.0$\%$\\
& $\#$4 & 6,497 & 6,502  & -0.1$\%$ & 11,101 &  2,501 & +77.5$\%$\\
& $\#$5 & 6,497 & 6,476  & +0.3$\%$ & 11,101 &  3,120 & +71.9$\%$\\
\hline
\end{tabular}
\label{table:comparison_objecitve_functions_values_test_problems}
\end{table}

Combining the results obtained for the different roads listed in  Table~\ref{table:comparison_objecitve_functions_values_test_problems}, we make an overall comparison between the two solvers. We observed that the HOPSPACK solver might yield a better or a worse solution than the solution obtained by the NOMAD solver. Thus, considering the tolerance of the difference in the optimum objective values obtained by the two solver, we count the number of times a solver wins with respect to the other solver. Table~\ref{table:overall_comparison_with_respect_to_objecitve_functions_values} shows the comparison of the solvers for different tolerance values of the difference in the optimum objective values. The $x \%$ tolerance of the difference in optimum objective values means if the optimum objective values obtained by the two solvers are in between $-x \%$ to $+x \%$, then the solvers yield the same solution (i.e., the two solvers tie), otherwise a positive percentage value indicates the HOPSPACK solver wins and a negative percentage value indicates the NOMAD solver wins.

\begin{table}[!ht]
\centering
\caption{Overall comparison of the HOPSPACK solver and the NOMAD solver with the optimum objective function values.}
\begin{tabular}{|>{\centering\arraybackslash} p{2.5 cm} | >{\centering\arraybackslash} p{2 cm}| >{\centering\arraybackslash} p{2 cm}| >{\centering\arraybackslash} p{2 cm}|}
\hline
Tolerance of the difference in the optimum costs & No. of times the NOMAD solver wins  & No. of times the HOPSPACK solver wins & No. of times the two solvers ties \\
\hline \hline
$\rpm $ 1$\%$ & 8 & 7  & 10\\
$\rpm $ 2$\%$   & 7 & 6  & 12\\
$\rpm $ 3$\%$  & 5 & 5  & 15\\
$\rpm $ 4$\%$  & 5 & 3  & 17\\
$\rpm $ 5$\%$  & 3 & 3  & 19\\
$\rpm $ 6$\%$  & 3 & 2  & 20\\
$\rpm $ 7$\%$  & 2 & 1  & 22\\
$\rpm $ 8$\%$  & 1 & 1  & 23\\
$\rpm $ 9$\%$  & 1 & 0  & 24\\
$\rpm $ 10$\%$  & 0 & 0  & 25\\
\hline
\end{tabular}
\label{table:overall_comparison_with_respect_to_objecitve_functions_values}
\end{table}

In Table~\ref{table:overall_comparison_with_respect_to_objecitve_functions_values}, we see that if the tolerance of difference in the optimum objective value is $\rpm 3 \%$ or above, then the two solver tie for more than $50\%$ test runs (i.e., more than 13 test runs among 25 test runs). We can also observe that for any case of the tolerance change in the optimum objective function value, the difference in the number of times the NOMAD solver wins and the number of times the HOPSPACK solver wins is at most 2. So in terms of optimum objective values obtained by the two solvers, the performance of both solvers are roughly equivalent. In Table~\ref{table:comparison_objecitve_functions_values_test_problems}, we can see that for all of the 25 test runs, the HOPSPACK solver required on average $78\%$ less black-box evaluation than the NOMAD solver. So the HOPSPACK solver is roughly five time faster than the NOMAD solver to compute the optimum solution.

\section{Conclusion}
Determining a good horizontal alignment is the first concern when designing a new road.  In this paper, we model horizontal alignment using the geometric specifications which are used by engineers in practice.  We further pursued a novel approach to address the horizontal alignment optimization problem. While most of the studies in the literature used heuristic based methods, we used a derivative-free optimization approach since it converges to a locally optimum solution.  Our method obtains a horizontal alignment that is locally optimal and a vertical alignment that is globally optimal. On our road samples, this locally optimal solution was on average $27\%$ cheaper than the solution provided by the engineers.  As the solution of the optimization model yields a practical horizontal alignment that satisfies geometric specifications and engineering requirements, the results represents significant savings in the final design of a new road.

Backward bends in a horizontal alignment (i.e., roads where the cheapest path from start to end requires moving away from the target destination for some portion of the horizontal alignment) are known to be a challenge to optimize \cite[page 123, Chapter 5]{Nicholson-1973}, \cite[page 21, Section 2.4.3]{JHA-06a}, \cite{Parker-1977}.  (Such behaviour is often referred to as backtracking.) However, our model can generate such alignments without difficult.

Although our proposed model works well for solving practical horizontal alignment optimization problems, it can be improved further for better precision and performance.
In our model formulation, we considered that the cross sections in a corridor are fixed and taken corresponding to the baseline alignment. However, cross sections should be always perpendicular to a horizontal alignment. When a horizontal alignment is significantly different from the baseline alignment, a set of new cross sections should be generated to increase the precision before calculating the corresponding vertical road profile.

A surrogate cost function is an approximation of the original cost function which is cheaper to compute. Using a surrogate cost function might reduce the solution time required to get a solution. (The NOMAD solver can exploit the usage of a non-adaptive surrogate.)

In our model, we only consider the construction costs to formulate our cost function. In the future, land acquisition costs could be incorporated by considering the unit cost of a piece of land corresponding to two consecutive cross sections in a corridor. We can also include pavement costs.

Both solvers need an initial starting point to start the algorithms. Since solutions obtained by both solvers are locally optimum, we can use multiple starting points (i.e.; multiple initial alignments) to obtain a better solution quickly. How to choose multiple alternative good alignments in a specified corridor is the subject of active investigation.

During the optimization process, at each iteration, both of the derivative free optimization solvers solve a large number of vertical alignment optimization problems (i.e., large scale mixed integer linear programming (MILP) problems). At the earlier stage of the optimization process (i.e., when the mesh size is coarse) we can relax some of the parameters of the vertical alignment optimization problem to get an approximate cost and then at the later stage (i.e., when the mesh size becomes relatively small) we can tighten the parameters to get accurate costs. This policy might reduce the solution time significantly. We can also use a warm start of the vertical alignment optimization problem when the horizontal alignments are close to each other to accelerate the vertical alignment optimization process. So the interconnection between the derivative free optimization solver and the MILP solver can be a potential way to reduce the solution time.

\section*{Acknowledgments}
This work was supported by the Natural Sciences and Engineering Research Council of Canada (NSERC) through a Collaborative Research and Development (CRD) grant sponsored by Softree Technical Systems Inc., a Discovery grant from the second author and a Discovery grant from the third author. Part of the research was performed in the Computer-Aided Convex Analysis (CA\textsuperscript{2}) laboratory funded by a Leaders Opportunity Fund (LOF) from the Canada Foundation for Innovation (CFI) and by a British Columbia Knowledge Development Fund (BCKDF).

The authors acknowledge the ideas, concepts, and implementation of the quasi network flow model for the earthwork problem with blocks, also known as the arc-flow model for the Earth Move problem with blocks, that are found in the work of Donovan Hare.

\appendix

\section{Vertical Alignment Optimization Model}\label{appendix}

In this appendix we give a very brief overview of the Quasi Network Flow (QNF) MILP used to solve the vertical alignment optimization subproblem.  The model is fully described in \cite{Hare-2014} and we refer there for details and definitions. The following equations are only provided to convince the reader that the model is indeed a MILP.

The variables are $f^k_{i,j,t}$, $V^+_i$, $V^-_i$, $u_i$, $a_{g,i}$, and binary variables $\nu_{i,l}$, $y_{k,t}$.
The objective function is
\begin{multline*}
\min \sum_{\substack{i \in \S \cup \B}} p_{i} V_{i}^{+}
+
\sum_{\substack{i \in \S \cup \W}} q_{i} V_{i}^{-}
+
\sum_{\substack{ i \in \S \\ t \in \T}} \left(c^r_{i,i-1} f^r_{i,i-1,t} + c^r_{i,i+1} f^r_{i,i+1,t} \right)\\
+
\sum_{\substack{j \in \B\\ t\in \T}} c\bar{d_j} \left( f_{j,\vartheta(j)-1,t}^b + f_{j,\vartheta(j)+1,t}^b \right) +
\sum_{\substack{j \in \W\\ t\in \T}} c\bar{d_j} \left( f_{\varphi(j)-1,j,t} + f_{\varphi(j)+1,j,t}^w \right)
.
\end{multline*}

The model includes numerous linear constraints like conservation of material
\begin{align*}
f^r_{i-1,i,t} + f_{i,i+1,t}^u + \sum_{\substack{j \in \B \\ \vartheta(j)=i}} f_{j,i+1,t}^b = f^r_{i,i+1,t} + f_{i-1,i,t}^l + \sum_{\substack{j \in \W \\ \vartheta(j)=i}} f_{i-1,j,t}^w, \\ \text{for all }t\in \T, i\in \S, \\
f^r_{i+1,i,t} + f_{i,i-1,t}^u + \sum_{\substack{j \in \B \\ \vartheta(j)=i}} f_{j,i-1,t}^b = f^r_{i,i-1,t} + f_{i+1,i,t}^l + \sum_{\substack{j \in \W \\ \vartheta(j)=i}} f_{i+1,j,t}^w, \\ \text{for all }t\in \T, i\in \S ,
\end{align*}
and volume balance constraints
\begin{align*}
&\sum_{t\in \T} f_{i,i-1,t}^u + f_{i,i+1,t}^u + f_{i,i,t}^u = V_i^+, &&\forall i \in \S,\\
&\sum_{t\in \T} f_{i-1,i,t}^l + f_{i+1,i,t}^l + f_{i,i,t}^l = V_i^-, &&\forall i \in \S,\\
&\sum_{t\in \T} f_{j,\vartheta(j)-1,t}^b + f_{j,\vartheta(j)+1,t}^b = V_j^+, &&\forall j \in \B,\\
&\sum_{t\in \T} f_{\varphi(j)-1,j,t} + f_{\varphi(j)+1,j,t}^w = V_j^-, &&\forall j \in \W.
\end{align*}

The volumes are related to the quadratic spline
\[
P(s) =  \begin{cases}
			P_1 (s) &\text{ if } s_{\d(1,1)} \leq s \leq s_{\d(1,n_1)},\\
			P_2 (s) &\text{ if } s_{\d(2,1)} \leq s \leq s_{\d(2,n_2)},\\
			\vdots & \\
			P_{\bar{g}}(s) & \text{ if } s_{\d(\bar{g},1)} \leq s \leq s_{\d(1,n_{\bar{g}})}.
		\end{cases}
\]
where $P_g(s) = a_{g,1}+a_{g,2}s+ a_{g,3} s^2$ by the offset value
\begin{align*}
&P(s_i) - h_i = u_i, &&\forall i \in \S.
\end{align*}

The smoothness of the quadratic spline is enforced by
\begin{align*}
P_{g-1}(s_{\d(g,1)}) &= P_g(s_{\d(g,1)}), &&\forall g \in \G \setminus \{1\}, &&\\
P'_{g-1}(s_{\d(g,1)}) &= P'_g(s_{\d(g,1)}) &&\forall g \in \G \setminus \{1\}, &&
\end{align*}
while grade constraints involves the derivative of the spline
\begin{align*}
G_L \leq P'_g(s_{\d(g,1)}) \leq G_U, &&\forall g \in \G \setminus \{1\}. &&
\end{align*}

Other constraints include fixed point constraints, side-slopes constraints, block constraints, and bound constraints. All these constraints are linear but may involve binary variables. We refer to \cite{Hare-2014} for their explicit definition and expression.

\bibliography{HAOpitmizationref}
\bibliographystyle{plain}

\end {document}